\documentclass[12pt]{article}
\usepackage[dvipsnames]{xcolor}
\usepackage{amsmath}
\usepackage{setspace}
\usepackage{bm}
\usepackage{bbm}
\usepackage{amssymb}
\usepackage{indentfirst}
\usepackage[superscript]{cite}
\usepackage{geometry}
\usepackage{mathtools}
\usepackage{tocloft}

\DeclareMathOperator{\sign}{sign}
\newcommand*{\citena}[1]{%
\begingroup
[\color{Green}
\romannumeral-`\x 
\setcitestyle{numbers}%
\cite{#1}%
\endgroup
]\ignorespacesafterend
}

\newcommand*{\eqrefe}[1]{%
\begingroup
(\color{BrickRed}
\romannumeral-`\x 
\setcitestyle{numbers}%
\ref{eq:#1}%
\endgroup
)\ignorespacesafterend
}

\newcommand*{\secrefe}[1]{%
\begingroup
(\color{Aquamarine}
\romannumeral-`\x 
\setcitestyle{numbers}%
\ref{#1}%
\endgroup
)\ignorespacesafterend
}

\newtagform{Tags}[\textcolor{BrickRed}]{\color{Black}(}{)}

\usepackage[super]{natbib}
\setcitestyle{super}

\newcommand{\ii}{\bm{i}}

\DeclareMathOperator{\csch}{csch}
\DeclarePairedDelimiter{\floor}{\lfloor}{\rfloor}

\begin{document}
\title{The Faulhaber Formula Analytic Continuation}
\date{January 12, 2021}
\author{Jose Risomar Sousa}
\maketitle
\usetagform{Tags}

\begin{abstract}
We extend the Faulhaber formula to the whole complex plane, obtaining an expression that fully resembles the Euler-Maclaurin summation formula, only it's exact. Thereafter, an expression for the generalized harmonic progressions valid in the whole complex plane is also derived. Lastly, we extend a formula for the Hurwitz zeta function valid at the negative integers, $\zeta(-k,b)$, to the whole complex plane, following a similar procedure.
\end{abstract}


\section{Introduction}
We present a new approach to the generalized harmonic numbers, $H_{-k}(n)$, based on the Faulhaber formula, which has the advantage of being valid in the whole complex plane, except a single point corresponding to the harmonic numbers ($k=-1$). The formula thus obtained allows us to deduce the exact value of the error term in the Euler-Maclaurin summation formula. Then, from this new formula, a new expression for $HP_{-k}(n)$ is derived as well.\\

Lastly, taking advantage of the same rationale, we also extend a formula for $\zeta(-k,b)$ only valid at the negative integers $-k$ to the whole complex plane. This analytic continuation bears a little resemblance to the Abel-Plana formula, but they are not the same, as evidenced by plots of their integrand functions.

\section{Generalized harmonic numbers} \label{Faulhaber}
It's possible to create a generalization of the Faulhaber formula and obtain $H_{-k}(n)$ for all complex $k$ except a single point. We start with the Faulhaber formula for odd powers:
\begin{equation} \nonumber
H_{-2k+1}(n)=\sum_{j=1}^{n}j^{2k-1}=\frac{n^{2k-1}}{2}+(2k-1)!\sum_{j=0}^{k-1}\frac{B_{2j}n^{2k-2j}}{(2j)!(2k-2j)!}
\end{equation}
\indent The idea is to use the analytic continuation of the Bernoulli numbers, achievable through the zeta function (as incredible as it may seem, we will see it works):
\begin{equation} \nonumber
\frac{B_{2j}}{(2j)!}=-2(-1)^j(2\pi)^{-2j}\zeta(2j) \text{}
\end{equation}\\
\indent And hence the formula becomes:
\begin{equation} \nonumber
H_{-2k+1}(n)=\frac{n^{2k}}{2k}+\frac{n^{2k-1}}{2}+2(2k-1)!(2\pi\ii)^{-2k}\zeta(2k)-2(2k-1)!\,n^{2k}\sum_{j=1}^{k}\frac{(2\pi\ii\,n)^{-2j}\zeta(2j)}{(2k-2j)!}
\end{equation}\\
\indent Now, let's recall two formulae we created for the zeta funcion out of the critical strip in \citena{AC}:
\begin{equation} \label{eq:zeta(k)pos}
\zeta(k)=\frac{1}{k!}\int_{0}^{1}\frac{\left(-\log{u}\right)^k}{(1-u)^2}\,du \text{, if } \Re{(k)}>1 \text{,}
\end{equation}

\begin{equation} \label{eq:zeta(k)neg}
\zeta(k)=-\frac{2(2\pi)^{k-1}}{k-1}\sin{\frac{k\,\pi}{2}}\int_{0}^{1}\frac{\left(-\log{u}\right)^{-k+1}}{(1-u)^2}\,du \text{, if } \Re{(k)}<0
\end{equation}\\
\indent Therefore, using equation \eqrefe{zeta(k)pos}, we have:
\begin{equation} \nonumber
\sum_{j=1}^{k}\frac{(2\pi\ii\,n)^{-2j}\zeta(2j)}{(2k-2j)!}=\int_{0}^{1}\sum_{j=1}^{k}\frac{(2\pi\ii\,n)^{-2j}}{(2k-2j)!}\frac{\left(-\log{u}\right)^{2j}}{(2j)!(1-u)^2}\,du
\end{equation}\\
\indent Replacing the sum inside the integral with a closed-form, we obtain its analytic continuation:
\begin{equation} \nonumber
\sum_{j=1}^{k}\frac{(2\pi\ii\,n)^{-2j}\zeta(2j)}{(2k-2j)!}=\frac{(2\pi\,n)^{-2k}}{2(2k)!}\int_{0}^{1}\frac{-2(2\pi n)^{2k}+\left(2\pi n-\ii\log{u}\right)^{2k}+\left(2\pi n+\ii\log{u}\right)^{2k}}{(1-u)^2}\,du
\end{equation}\\
\indent This integral can be generalized for any parity (it's used later down the line in section \secrefe{Hur_AC}) and transformed in a few ways, so as to rid it of the non-real parts:
\begin{multline} \label{eq:soma_zetas_pares}
\sum_{j=1}^{\floor{k/2}}\frac{(2\pi\ii\,n)^{-2j}\zeta(2j)}{(k-2j)!}=\frac{(2\pi\,n)^{-k}}{2\,k!}\int_{0}^{1}\frac{-2(2\pi n)^{k}+\left(2\pi n-\ii\log{u}\right)^{k}+\left(2\pi n+\ii\log{u}\right)^{k}}{(1-u)^2}\,du\\=-\frac{\sign{\Re{(n)}}\,\pi\,n}{2\,k!}\int_{0}^{\pi/2}\left(\sec{v}\csch{(\pi\,n\tan{v})}\right)^2\left(1-\frac{\cos{k\,v}}{(\cos{v})^{k}}\right)\,dv\\=-\frac{\sign{\Re{(n)}}\,\pi\,n}{2\,k!}\int_{0}^{\infty}\left(\csch{(\pi\,n\,v)}\right)^2\left(1-\left(1+v^2\right)^{k/2}\cos{\left(k\arctan{v}\right)}\right)\,dv
\end{multline}\\
\indent The first integral is always equal to the two others, if $\Re{(n)}\ge 0$, while for some values of $k$ and $n$ the first integral may diverge, whereas the two last may converge (e.g., non-real $k$ and $n$, with $\Re{(k)}<0$ and $\Re{(n)}<0$).\\

\indent Therefore, since $n$ is a positive integer, the analytic continuation of the Faulhaber formula, which, surprisingly, holds in the whole complex plane (except $k=0$) is:
\begin{multline} \nonumber
H_{-2k+1}(n)=\frac{n^{2k}}{2k}+\frac{n^{2k-1}}{2}+2(2k-1)!(2\pi)^{-2k}(\cos{k\pi})\zeta(2k)\\+\frac{\pi\,n^{2k+1}}{2k}\int_{0}^{\pi/2}\left(\sec{v}\csch{(\pi\,n\tan{v})}\right)^2\left(1-\frac{\cos{2k\,v}}{(\cos{v})^{2k}}\right)\,dv \text{,}
\end{multline}\\
\noindent where $\cos{k\pi}$ is the real part of $\ii^{-2k}$, a transformation that's necessary for the formula to be right.\\

Finally, since the formula now holds for any $k$, we can transform it into a proper form, by making $2k-1$ equal to $k$. This generalized Faulhaber formula holds for all complex $k$, except $-1$:
\begin{equation} \label{eq:HN_final}
\sum_{j=1}^{n}j^{k}=\frac{n^{k+1}}{k+1}+\frac{n^{k}}{2}+\zeta(-k)+\frac{\pi\,n^{k+2}}{k+1}\int_{0}^{\pi/2}\left(\sec{v}\csch{\left(\pi\,n\tan{v}\right)}\right)^2\left(1-\frac{\cos{(k+1)\,v}}{(\cos{v})^{k+1}}\right)\,dv \text{,}
\end{equation}\\
\noindent where $\zeta(-k)$ came from:
\begin{equation} \nonumber
2\,k!(2\pi)^{-(k+1)}\cos{\frac{(k+1)\pi}{2}}\zeta(k+1)=\zeta(-k) \text{,}
\end{equation}
\noindent which is nothing but the Riemann functional equation.\\

Now it's easy to see what's being discarded when one analytically continues the Riemann zeta function from $\Re{(k)}>1$ to the whole complex plane.\\

\indent An alternative representation, which can be more useful at times, is shown below:
\begin{equation} \nonumber
\sum_{j=1}^{n}j^{k}=\frac{n^{k+1}}{k+1}+\frac{n^{k}}{2}+\zeta(-k)-\frac{2\pi\,n^{k+2}}{k+1}\int_{0}^{\infty}\left(-2+(1+\ii\,x)^{k+1}+(1-\ii\,x)^{k+1}\right)\frac{e^{-2\pi\,n\,x}}{\left(1-e^{-2\pi\,n\,x}\right)^2}\,dx 
\end{equation}

\section{Generalized harmonic progressions}
Through the generalized Faulhaber formula, we can also obtain an expression for $HP_{-k}(n)$ that holds for all complex $k$, except the harmonic progression (of order 1).
\begin{equation} \nonumber
HP_{-2k+1}(n)=\sum_{i=1}^{n}(i+b)^{2k-1}=(2k-1)!\,b^{2k-1}\sum_{j=0}^{2k-1}\frac{b^{-j}H_{-j}(n)}{j!(2k-1-j)!}
\end{equation}\\
\indent Therefore, we replace $H_{-j}(n)$ with its equivalent formula from section \secrefe{Faulhaber}, equation \eqrefe{HN_final} plus a transformation similar to \eqrefe{soma_zetas_pares}, and then need to resolve the below expression and find closed-forms for each one of its parts:
\begin{multline} \nonumber 
(2k-1)!\,b^{2k-1}\sum_{j=0}^{2k-1}\frac{b^{-j}}{j!(2k-1-j)!}\left(\frac{n^{j+1}}{j+1}+\frac{j^{k}}{2}+\zeta(-j)\right)\\-(2k-1)!\,b^{2k-1}\sum_{j=0}^{2k-1}\frac{b^{-j}(2\pi)^{-j-1}}{(j+1)!(2k-1-j)!}\int_{0}^{1}\frac{-2(2\pi n)^{j+1}+\left(2\pi n-\ii\log{u}\right)^{j+1}+\left(2\pi n+\ii\log{u}\right)^{j+1}}{(1-u)^2}\,du 
\end{multline}\\
\indent Only one of these parts is difficult to figure, one is trivial and the rest falls into the same somewhat simple pattern. First, since the zeta function vanishes at the even negative integers, we have:
\begin{multline} \nonumber 
(2k-1)!\,b^{2k-1}\sum_{j=0}^{2k-1}\frac{b^{-j}\zeta(-j)}{j!(2k-1-j)!}=(2k-1)!\,b^{2k-1}\left(-\frac{1}{2(2k-1)!}+\sum_{j=1}^{k}\frac{b^{-2j-1}\zeta(-2j+1)}{(2j-1)!(2k-2j)!}\right)
\end{multline}\\
\indent Replacing the zeta function at the negative integers with its expression from formula \eqrefe{zeta(k)neg}:
\begin{equation} \nonumber 
\sum_{j=1}^{k}\frac{b^{-2j-1}\zeta(-2j+1)}{(2j-1)!(2k-2j)!}=\int_{0}^{1}\sum_{j=1}^{k}\frac{b^{-2j-1}}{(2j-1)!(2k-2j)!}\frac{(-1)^j(2\pi)^{-2j}}{j}\frac{(\log{u})^{2j}}{(1-u)^2}\,du \text{,}
\end{equation}\\
\noindent and therefore:
\begin{multline} \nonumber 
(2k-1)!\,b^{2k-1}\sum_{j=0}^{2k-1}\frac{b^{-j}\zeta(-j)}{j!(2k-1-j)!}=-\frac{b^{2k-1}}{2}\\+\frac{b^{2k}}{k}\int_{0}^{1}\left(-1+\left(1+\left(\frac{\log{u}}{2\pi\,b}\right)^2\right)^k\cos{\left(2\,k\arctan{\frac{\log{u}}{2\pi\,b}}\right)}\right)\frac{1}{(1-u)^2}\,du 
\end{multline}\\
\indent If $\Re{(b)}>0$, we can simplify this integral a little:
\begin{multline} \nonumber 
(2k-1)!\,b^{2k-1}\sum_{j=0}^{2k-1}\frac{b^{-j}\zeta(-j)}{j!(2k-1-j)!}=-\frac{b^{2k-1}}{2}-\frac{\pi\,b^{2k+1}}{2\,k}\int_{0}^{\pi/2}\left(\sec{v}\csch{\left(\pi\,b\tan{v}\right)}\right)^2\left(1-\frac{\cos{2\,k\,v}}{(\cos{v})^{2\,k}}\right)\,dv
\end{multline}\\
\indent Going back to the original expression, the integral part can be simplified as follows:
\begin{multline} \nonumber 
-(2k-1)!\,b^{2k-1}\sum_{j=0}^{2k-1}\frac{b^{-j}(2\pi)^{-j-1}}{(j+1)!(2k-1-j)!}\int_{0}^{1}\frac{-2(2\pi n)^{j+1}+\left(2\pi n-\ii\log{u}\right)^{j+1}+\left(2\pi n+\ii\log{u}\right)^{j+1}}{(1-u)^2}\,du=\\
-\frac{(2\pi)^{-2k}}{2k}\int_{0}^{1}\frac{-2(2\pi(n+b))^{2k}+\left(2\pi(n+b)-\ii\log{u}\right)^{2k}+\left(2\pi(n+b)+\ii\log{u}\right)^{2k}}{(1-u)^2}\,du \text{,}
\end{multline}\\
\noindent which for $\Re{(b)}>0$ can be turned into a simpler integral.\\

Though using the original integrals the final formula may hold for many other values of $b$, the below expression holds for all complex $k$ (except $-1$), and all $b$ with $\Re{(b)}>0$. Introducing the sign of $\Re{(b)}$ in the integral per relation \eqrefe{soma_zetas_pares} makes it more accurate, but the equation still doesn't hold always when $\Re{(b)}<0$:
\begin{multline} \nonumber 
\sum_{j=1}^{n}(j+b)^k=\frac{(n+b)^{k+1}-b^{k+1}}{k+1}+\frac{(n+b)^k-b^k}{2}\\+\frac{\sign{\Re{(b)}}\,\pi}{k+1}\int_{0}^{\pi/2}(\sec{v})^2\left((n+b)^{k+2}\csch^2{\left(\pi(n+b)\tan{v}\right)}-b^{k+2}\csch^2{\left(\pi\,b\tan{v}\right)}\right)\left(1-\frac{\cos{(k+1)\,v}}{(\cos{v})^{k+1}}\right)\,dv
\end{multline}\\
\indent Now, more strikingly, we can transform this formula further using equation \eqrefe{Hurwitz_int}. It should hold always when $k \neq -1$ and $\Re{(b)}>0$:
\begin{multline}  
\sum_{j=0}^{n}(j+b)^k=\frac{(n+b)^{k+1}}{k+1}+\frac{(n+b)^k}{2}+\zeta(-k,b)\\+\frac{\sign{\Re{(b)}}\,\pi(n+b)^{k+2}}{k+1}\int_{0}^{\pi/2}\left(\sec{v}\csch{\left(\pi(n+b)\tan{v}\right)}\right)^2\left(1-\frac{\cos{(k+1)\,v}}{(\cos{v})^{k+1}}\right)\,dv
\end{multline}

\section{Hurwitz zeta analytic continuation} \label{Hur_AC}

Based on a relation between $\zeta(-k,b)$ and $\mathrm{Li}_{k}(e^{\pm 2\pi\ii\,b})$, available from the literature, a formula for $\zeta(-k,b)$ that holds at the non-positive integers $-k$ was created in \citena{Lerch}, which holds for every complex $b$:
\begin{equation} \label{eq:Hur_neg_int}
\zeta(-k,b)=-\frac{b^{k+1}}{k+1}+\frac{b^k}{2}+2\,k!\,b^{k+1}\sum_{j=1}^{\floor{(k+1)/2}}\frac{(-1)^{j}(2\pi\,b)^{-2j}\zeta(2j)}{(k+1-2j)!}
\end{equation}\\
\indent So, now we can easily create the analytic continuation of the Hurwitz zeta function using the formula \eqrefe{Hur_neg_int} and the analytic continuation from equation \eqrefe{soma_zetas_pares}.\\

The below integral representation for the Hurwitz zeta function holds for every complex $k\neq 1$ and $b$ such that $\Re{(b)}>0$ (it doesn't hold always when $\Re{(b)}<0$):
\begingroup
\small
\begin{equation} \label{eq:Hurwitz_int} 
\zeta(k,b)=\frac{b^{-k+1}}{k-1}+\frac{b^{-k}}{2}+\frac{\sign{\Re{(b)}}\,\pi\,b^{-k+2}}{k-1}\int_{0}^{\pi/2}\left(\sec{v}\csch{\left(\pi\,b\tan{v}\right)}\right)^2\left(1-\frac{\cos{(k-1)\,v}}{(\cos{v})^{-k+1}}\right)\,dv
\end{equation}
\endgroup\\
\indent Finally, using equation \eqrefe{Hurwitz_int}, an expression for the zeta function that holds always, except at $1$, is:
\begin{equation} \nonumber
\zeta(k)=\frac{1}{k-1}+\frac{1}{2}+\frac{\pi}{k-1}\int_{0}^{\pi/2}\left(\sec{v}\csch{\left(\pi\tan{v}\right)}\right)^2\left(1-\frac{\cos{(k-1)\,v}}{(\cos{v})^{-k+1}}\right)\,dv
\end{equation}


\end{document}